\documentclass[12pt]{article}
\usepackage{amsmath,amsfonts,amssymb,amsthm}

\def\noi{\noindent}

\def\pf{\noi{\bf Proof.\ \,}}
\def\eop{{$\square$}}

\def\labtt#1{\label {#1}}
\def\labttr#1{\label {#1}\rm }

\def\g{\gamma}
\def\d{\delta}

\def\o{\omega}

\def\vep{\varepsilon}

\def\FF{{\mathbb F}}

\def\QQ{{\mathbb Q}}
\def\RR{{\mathbb R}}

\def\ZZ{{\mathbb Z}}

\def\la{\langle}
\def\ra{\rangle}
\def\<{\langle}
\def\>{\rangle}

\def\half{{\frac 12}}
\def\fourth{{\frac 14}}

\def\elm#1#2{2^{#1-1}-2^{#1-#2-1}}
\def\elp#1#2{2^{#1-1}+2^{#1-#2-1}}

\def\dual#1{#1^*}        

\def\kron#1#2{\delta_{#1#2}}  

\def\lpt{L^+(t)}
\def\lmt{L^-(t)}
\def\lvept{L^\vep (t)}

\def\mv#1{MinVec(#1)}

\def\explus#1{2^{1+2#1}_+}

\def\ratholoex#1{2^{1+2#1}_+\Omega^+(2#1,2)}

\def\dg#1{{\cal D}({#1})}  

\def\bw#1{BW_{{2^{#1}}}}

\def\brw#1{BRW^+(2^{#1})}

\def\gd#1{G_{2^{#1}}} 
\def\rd#1{R_{2^{#1}}}

\def\pow#1{{\mathcal {P}}(#1)}  

\begin{document}

\newtheorem{thm}{Theorem}[section]
\newtheorem{prop}[thm]{Proposition}
\newtheorem{lem}[thm]{Lemma}
\newtheorem{rem}[thm]{Remark}
\newtheorem{coro}[thm]{Corollary}
\newtheorem{conj}[thm]{Conjecture}
\newtheorem{de}[thm]{Definition}
\newtheorem{hyp}[thm]{Hypothesis}

\newtheorem{nota}[thm]{Notation}
\newtheorem{ex}[thm]{Example}
\newtheorem{proc}[thm]{Procedure}  
\def\refp#1{(\ref {#1})}
\def\refpp#1{(\ref {#1})}

\begin{center}
{\Large \bf  Midwest cousins of Barnes-Wall lattices}
\medskip 
\vskip .5cm

{filename: {\bf  mwfirstcousinsrevision6sep09}; 

version 6 September, 2009; 

this replaces previous versions of  ``Midwest cousins of Barnes-Wall lattices''  }

{\Large    }

\vspace{10mm}
Robert L.~Griess Jr.
\\[0pt]
Department of Mathematics\\[0pt] University of Michigan\\[0pt]
Ann Arbor, MI 48109 USA  \\[0pt]
\vskip 1cm 

\end{center}

\newpage

\begin{abstract} 
Given a rational lattice and suitable set of linear transformations, 
we construct a cousin lattice.   
Sufficient conditions are given for integrality, evenness and unimodularity. 
When  the input is a Barnes-Wall lattice, we get multi-parameter series of cousins.  There is a subseries consisting of unimodular  lattices which have  ranks $2^{d-1}\pm 2^{d-k-1}$, for odd integers $d\ge 3$ and integers $k=1,2, \cdots , \frac {d-1}2$.  
Their  minimum norms are  moderately high:  $2^{\lfloor \frac d2 \rfloor -1}$.   
\end{abstract}

\vskip 5cm 


\newpage 

\tableofcontents

\section{Introduction} 

In this article, lattice means a finite rank free abelian group with rational-valued positive definite symmetric bilinear form.

We develop a general lattice construction method which is inspired by finite  group theory.  
We call it a {\it midwest procedure} 
because 
many significant developments in finite group theory took place in the American midwest during the late twentieth century, especially in Illinois, Indiana,  Michigan, Ohio and Wisconsin.  

The idea is to start with a lattice $L$ and take a finite subgroup $F$ of $O(\QQ\otimes L)$.  
In the rational span of $F$ in $End(\QQ \otimes L)$, 
we take an element  $h$.  
We define a new lattice, $L'$, in some way using $L$ and $h$, 
for example $L \cap Ker(h)$, $\dual L \cap Ker(h)$, $Lh$, \dots  ,  or sums of such things.  
After finitely many repetitions of this procedure, 
the sequence $L, L', \dots $ arrives at a new lattice, which is called a {\it midwest cousin} of $L$.   In this article, we restrict this procedure to the midwest cousins defined in \refpp{mc}.

In \refpp{mcbw}, we specialize further to the dimension $2^d$ Barnes-Wall lattices $\bw d$ and the Bolt-Room-Wall groups $\brw d$, of shape $\ratholoex d$, which are the full isometry groups of $\bw d$ if $d\ne 3$.    The sophisticated groups $\brw d$ help us manage the linear algebra and combinatorics.  
We create multi-parameter series of cousin lattices, called {\it the first cousins of the Barnes-Wall lattices}.  The dimension of a first cousin is $2^{d-1}\pm 2^{d-k-1}$, for some $k\in \{1, 2, \cdots , \lfloor \frac d2 \rfloor \}$.  
The auxiliary finite isometry groups $F_i$ are cyclic groups of orders 2 and 4.  
When $d$ is odd and $d-2k\ge 3$, 
the minimum norms are  $2^{\lfloor \frac d2 \rfloor -1}$ and the lattices are even and unimodular.  We include a partial analysis of minimal vectors.

We are grateful to the University of Michigan,  National Cheng Kung University, Zhejiang University,  and the U. S. National Science Foundation  for financial support (NSF DMS-0600854).  We thank Harold N. Ward for useful discussions and the referee for many helpful comments.

\subsection{Conventions and List of Notations}

Group elements and endomorphisms usually act on the right.  
Table 1 summarizes notations.  
 An appendix to this article summarizes background.  For more details, see \cite{poe,bwy,ibw1}. 
The upcoming book \cite{grgl} may be helpful.

\begin{table}[htdp]
\caption{List of Notations, Part 1}
\begin{center}
\begin{tabular}{|c|c|c|}
\hline 
Notation & Summary &  Comments \cr 
\hline 
\hline
$\brw d$ & the Bolt-Room-Wall group, $\ratholoex d$ & \cr
\hline 
$\bw d$ &the Barnes-Wall lattice of rank $2^d$ & \cr 
\hline 
BW-level& & \refpp{levels} \cr
\hline
commutator density && \refpp{commdens}\cr 
\hline
$Mod(D,-)$ &category of modules for $D \cong 
Dih_8$& \cr 
& where central involution acts as $-1$& \cr
\hline  
core & $S_1\cap \cdots \cap S_r$ as in cubi sum (below)& \refpp{core} \cr 
\hline 
cubi sum & $S_1+\cdots + S_r$, $S_i$ affine  codimension 2 &\refpp{core} \cr
& subspaces in $\FF_2^d$  so that & \cr
& $codim(S_1\cap \cdots \cap S_r)=2r$ & cubi theory \cite{ibw1} \cr
\hline 
defect & invariant of an involution in $\brw d$ & \refpp{traces}\refpp{rem1},\refpp{eq1} \cr
\hline 
$\vep_S$ & $v_i \mapsto v_i, -v_i$, as $i \not \in S, i \in S$& \cr
\hline 
fourvolution & an isometry of order 4 whose square is $-1$ & \cr 
\hline 
frame, lower frame& & \refpp{lowerframe}\cr
\hline 
$G, \gd d$ & $\brw d$, a subgroup of $O(\bw d )$ &\refpp{bw0} \cr
\hline
Jordan number, JNo & & \refpp{jno}\cr
\hline 
$k^{th}$ layer & $L(k)/L(k-1)$ & \cr 
\hline
level & least $\ell$ so that $2^\ell x$ has integer coordinates & \refpp{levels} \cr 
\hline 
level sublattice & & \refpp{layers}\cr
\hline 
$\lvept$ & eigenlattice for involution $t$ & \cr
\hline
$k^{th}$ level, $L(k)$& the set of lattice elements of level at most $k$ & \cr
\hline 
long codeword& $RM(2,d)$ codeword of weight more than $2^{d-1}$ & \cr 
\hline $\lpt, \lmt$&eigenlattices for involution $t$&\cr
\hline 
lower element & element of $\gd d$ contained in $\rd d$& \cr
\hline
lower frame & &\refpp{lowerframe} \cr
\hline
$MC(L,t,f,\vep)$ &a cousin lattice  &\refpp{mc} \cr
\hline 
$MC(\bw d, t, f, \vep)$ &a cousin lattice  & \refpp{mcbw}  \cr
\hline 
$MC_1(d,k,\vep)$ &a cousin lattice  & \refpp{mcbw} \cr
\hline 
$\mu (L)$ & the minimum norm in the lattice $L$ & \cr 
\hline
$O(L)$ & isometry group of quadratic space $L$ & \cr 
\hline 
\end{tabular}
\end{center}
\label{notationspart1}
\end{table}%


\begin{table}[htdp]
\caption{List of Notations, Part 2}
\begin{center}
\begin{tabular}{|c|c|c|}
\hline 
$O_p(X)$& the largest normal $p$-subgroup of &\cr
&the group $X$ ($p$ prime)& \cr 
\hline
$O_{p'}(X)$& the largest normal subgroup of &\cr
& the group $X$ of order prime to $p$ & \cr 
\hline
$\pow X$ & the power set of the set $X$ & \refpp{levels}, \refpp{tau-1} \cr 
\hline 
quotient code &quotient space of a code which has code structure & \refpp{quotientcode}\cr
\hline 
$R, \rd d$ &$O_2(\gd d)$&\refpp{bw0} \cr
\hline 
$RM(k,d)$ & the Reed-Muller code of length $2^d$ & \refpp{rmsetup} \cr 
\hline 
RM-level & & \refpp{levels} \cr
\hline
$sBW$, $ssBW$ & scaled, suitably scaled BW lattice & \cite{bwy} \cr
\hline
short codeword, &codeword in $\FF_2^n$ of weight $< \half n$ & \cr 
short involution & & \refp{termsrm2d} \cr
\hline
split, nonsplit  &involution of $\gd d$ which&\cr
& centralizes, does not centralize, &\cr
&a lower elementary abelian $2^{d+1}$&\refpp{ibw1review}\cr 
\hline 
standard frame, basis&&\refpp{lowerframe}\cr
\hline
standard generators & certain set of $2^{-m}v_A$ in $\bw d$ & \refpp{standardgenset} \cr 
\hline 
$t=\vep_A$& a diagonal involution in $\brw d$& \cr
\hline 
$\tau_\o , \o \in \Omega$ & translation by $\o$ on  $\Omega$ or $V:=\QQ \otimes \bw d$ &\refpp{bw0} \cr 
\hline 
$\tau (core(Z))$ & the group $\{\tau_c \mid c\in core(Z)\}$ & \cr
\hline   
$Tel(L,E), E$ abelian & total eigenlattice on lattice $L$, the sum of eigenlattices & \cr 
\hline 
$Tel(L,t), t$ involution& total eigenlattice on lattice $L$, $\lpt \perp \lmt$ & 
\refpp{jnok=0}, \refpp{mc2} \cr 
\hline 
$top(x)$ & part of vector $x$ representing the & \cr
& highest power of 2 in denominator
 & \refpp{top}  \cr 
 \hline
 top closure & $top(x)$ is in lattice if $x$ is in lattice& \refpp{top}, \refpp{topnotin}\cr 
\hline
upper element & element of $\gd d$ not contained in $\rd d$& \cr
\hline 
$V^\vep (t)$, $V$ & $\QQ \otimes \lvept$, $V:=\QQ \otimes L$ &\refpp{mc}  \cr 
\hline 
$v_i, v_X \in \RR^{\Omega}$ & $(v_i,v_j)=2^{\lfloor \frac d2 \rfloor} \kron ij$; $v_X:=\sum_{i \in X}v_i$&\refpp{standardgenset} \cr 
\hline 
$Z, Z+\Omega \in RM(2,d)$&  weight $2^{d-1}\pm 2^{d-k-1}$  codewords &
\refpp{rem1}, \refpp{core} \cr
\hline 
2/4, 3/4 generation & a property of some objects in $Mod(D,-)$& \refpp{2/4}\cr 
\hline 
$\Omega, \Omega_d$& index set for orthogonal basis of $\RR^{2^d}$ &\refpp{bw0}  \cr 
\hline 
\end{tabular}
\end{center}
\label{notationspart2}
\end{table}%

\section{Involutions on Barnes-Wall lattices}  

We use the notations and results of \cite{bwy} and 
\cite{ibw1}, which are recommended for background.  

\begin{de}\labttr{ibw1review} 
We recall   that an involution in $\brw d$ 
has trace 0 if and only if it is conjugate to its negative in $\brw d$ (equivalent, conjugate to its negative by an element of $\rd d$ 
\cite{grmontreal,bwy,ibw1}).  

 An involution in $\brw d$ is {\it split} if it centralizes a maximal elementary abelian subgroup of $\rd d$ and is {\it nonsplit} otherwise.  
\end{de}

For a summary of properties and classification of such involutions, see \cite{ibw1}  Appendix: 
About BRW groups.  We have changed some terminology since that article.   We mention one often-used result.

\begin{thm}\labtt{traces}(i) 
If $g\in \brw d$, then the trace of $g$ on the natural $2^d$-dimensional module is 0 or is $\pm 2^e$ if $g$ has nonzero trace, where $2e$ is the dimension of the fixed point subspace for the conjugation action of $g$ on $\rd d/Z(\rd d)$.  

(ii)  Suppose that $g\in \brw d$ is an involution.  
The defect $k$ of $g$ satisfies $e+k=d$.  
The multiplicities of eigenvalues $\pm 1$ are (up to transposition) $\elp dk, \elm dk$, respectively.  
\end{thm}

\begin{rem} 
\labttr{rem1} 
Let $A\in RM(2,d)$ be a short codeword of defect $k$ \refpp{core}.   
Throughout this article, we shall work with involutions of the form $t:=\vep_A$.  Its trace is $2^{d-k}$.  Let $A=A_1+\cdots +A_k$ be a cubi sum \refpp{core}.  
The affine subspace $core(A)=core(Z)=\cap_i A_i$ is $(d-2k)$-dimensional.  
For $c\in \Omega$, the  corresponding translation map is $\tau_c$.  
If $c\in core(A)$, we call $\tau_c$ a {\it core translation}, so when $core(A)$ contains the origin, we get a group of translations.  
Let $\tau_c$ be a nonidentity core translation.  
Observe that if we take any hyperplane $H$ which contains no translate of $c$, then  $f:=\vep_H\tau_{c}$ is a fourvolution which commutes with $t$.  
\end{rem}

\subsection{Involutions on Barnes-Wall lattices mod 2: JNo}  

We begin by studying the Jordan canonical form of involutions on the Barnes-Wall lattice modulo 2.  We derive applications to discriminant groups and lattice constructions.  

\begin{de}\labttr{jno} The {\it Jordan number} of an involution acting on a finite rank abelian group $A$ is the number of degree 2 Jordan blocks in its canonical form on $A/2A$.  We write $JNo(t)$ or $JNo(t,A)$ for the Jordan number of $t$.  
\end{de}

\begin{lem}\labtt{jnok=0} On $\bw d$, the Jordan number for $-1$ is 0 and 
the Jordan number is $2^{d-2}$ for a lower noncentral involution.  
\end{lem}
\pf
The first statement is obvious.  The second follows since $|\bw d :Tel(t)|=2^{2^{d-2}}$ for lower involutions $t$.  
See \cite{bwy}.  
\eop

\begin{nota}\labttr{jnodk} 
In this section, the notations of \refpp{jno} will stand for  lattices 
(which often  are sBWs) and the involutions will be isometries of them.   
Let $L$ be a sBW lattice of rank $2^d$.  
If $t\in O(L)$ is an involution, as before, we let $JNo(t)$ be its Jordan number \refpp{jno}.   
Because of \refpp{jnok=0}, we assume that the defect $k$ is positive, i.e., that the involution is upper.  
If $2k<d$, there exists a lower dihedral group in $C_{\gd d}(t)$.  
\end{nota}

Theorem \refpp{jno=elm} is  the main goal of this section.

\begin{lem}\labtt{nonsplitfull} 
If $t$ is a nonsplit involution, it has full Jordan number, i.e., $JNo(t)=2^{d-1}$.  
\end{lem} 
\pf  A nonsplit involution is upper.  
By \cite{ibw1}, there exists a lower dihedral group $D$ so that $t$ normalizes $D$ and effects an outer automorphism on $D$, say by transposing a set of generators $u, v$.  
Using 2/4 generation of $L$ with respect to $D$, we get $L=L^+(u)\oplus L^+(v)$ for a generating pair of involutions $u, v$ so that $u^t=v$.  Then obviously $L$ is a free $\ZZ \la t\ra$-module, so we are done.  
\eop

\begin{lem}\labtt{eq1}
If $t$ centralizes a lower dihedral group, $JNo(t)=JNo(t')+JNo(t'')$, where $t', t''$ are defect $k$ involutions on  sBW lattices of rank $2^{d-1}$.  
\end{lem}
\pf 
We may choose such a lower dihedral group $D$ to satisfy $D\cap [R,t]=Z(R)$.  
Use the 2/4 property to get that $t$ preserves each direct summand in  $L=L^+(u)\oplus L^+(v)$ for a generating pair of involutions $u, v$ of $D$ (the summands are sBW).  In the notation of \cite{bwy}, there exists 
a group $Q\cong \explus {(d-1)}$ in $\brw d$ which acts trivially on $L^-(u)$ and as a lower group on $L^+(u)$.  Since the action of $t$ on $R$ has defect $k$, the action of $t$ on $Q$ has defect $k$.  We may therefore apply induction to the restriction of $t$ to the summand $L^+(u)$.  A similar argument applies to $L^+(v)$.    
\eop

\begin{lem}\labtt{d=2} When $(d,k)=(2,1)$ and $t$ is an upper  involution, 
$JNo(t)=1$ when $t$ has nonzero trace and $JNo(t)=2$ when $t$ has trace zero.  
\end{lem}
\pf 
We refer to \cite{bwy} for a discussion of involutions in $\brw 2 \cong W_{F_4}$.

Suppose that the involution has nonzero trace.   
Since its trace is $\pm 2$, we may assume that it is 2, whence $t$ is a reflection.  Then the statement is obvious since reflections induce transvections on the lattice mod 2. 

For $d=2$, 
if an involution  is upper and  nonsplit, we may quote 
\refpp{nonsplitfull}.  
For $d=2$, 
if an involution  is upper and  split, it has nonzero trace and we may quote the previous paragraph.  
\eop

\begin{lem}\labtt{ineq1} If $t$ has nonzero trace, $JNo(t) \le \elm d k$.  
\end{lem}
\pf We may assume that $tr(t)>0$.  
Let $h$ be the dimension of fixed points for $t$ on $L/2L$.  Then $h+JNo(t)=2^d$.   Since the 1-eigenlattice for $t$ has rank $\elp d k$ and is a direct summand of $L$, we have $h\ge \elp d k$.
\eop

\begin{lem}\labtt{ineq2}  Suppose that the upper involution $t$ lies in a subgroup $S$ of $G$  of order $2n$, $n$ odd, and that every nonidentity element of $S$ of order dividing $n$ 
has the same fixed point subspace, of dimension $2e$, on $R/R'$.  Assume further that $t$ inverts a nonidentity odd order element of $S$. 
Then $JNo(t) \ge 2^{d-1} - \half (\frac {2^d-2^e} {n} + 2^e )$.
\end{lem}
\pf
Such a group $S$ has a normal subgroup of order $n$.  Call it $C$.  
Then every nonidentity element of $C$ has trace $\pm 2^e$ on $L$ \refpp{traces}.  It follows that the eigenlattice $M$ of $C$-fixed points has rank $\frac 1 n (2^d +(n-1)2^e)=\frac 1 n (2^d-2^e +n2^e)$.  
On the annihilator $N:=L \cap M^\perp$, $C$ acts faithfully on every constituent, 
and since $t$ inverts a nonidentity element of $C$, $N/2N$ is a free 
$\la t \ra$-module, whence $JNo(t) \ge \half rank(N)= \half (2^d-rank(M))$.  
\eop

Next, we deal with the situation when $t$ does not centralize a lower dihedral group.

\begin{lem}\labttr{ineq3}  We use the hypotheses and notation of \refpp{ineq2}.  

(i) Suppose that $d$ is even,  $n=2^{\frac d 2}+1$ and $e=0$.    Then 
$JNo(t)\ge   2^{d-1} - 2^{\frac d 2 -1}$.  

(ii) Suppose that $d$ is odd, $n=2^{\frac {d-1} 2}+1$ and $e=1$.    Then 
$JNo(t)\ge   2^{d-1} - 2^{\frac {d-1} 2}$.  
\end{lem}
\pf
Straightforward with \refpp{ineq2}.  
\eop

\begin{lem}\labtt{singer1} Suppose that $m\ge 1$, $2r\ge 4m \ge 4$ and that  $u$ is an involution in 
$\Omega ^+(2r,2)$ with commutator submodule of dimension $2m$ 
on its natural module $W:=\FF_2^{4m}$.  Assume that $W(u-1)$  is a totally singular subspace.  
Let $n=2^{2m}-1$.  

Then $u$ is in a group $P$ of order $2n$, where $P$ contains a Singer cycle  $C$ 
in a natural $GL(2m,2)$-subgroup of $\Omega ^+(2r,2)$ (so $C$ is a normal subgroup of $P$).  Also $P$ has the property that the nonidentity elements of $C$ have the same fixed point subspace on $\FF_2^{2r}$. 
\end{lem}
\pf 
Recall properties of the normalizer of a Singer cycle in classical groups, \cite{Hup}.  Without loss, we may assume that $2r=4m$.  

Suppose that we are given a pair of maximal totally singular subspaces, 
$W_1, W_2$ in $W$ 
 such that $W=W_1\oplus W_2$.   Let $H$ be the common stabilizer of $W_1$ and $W_2$.  So, $H\cong GL(2m,2)$.  
Let $P$ be the subgroup of the normalizer of a Singer cycle in 
$H$ corresponding to the Singer cycle and the group of field automorphisms of order 2.  It has order $2n$ and its involutions invert nonidentity elements of $C$ so have Jordan number $2m$ on $W$.  
If $u$ is conjugate to such an involution, we are done.  
There are two conjugacy classes of involutions in $\Omega ^+(2m,2)$ with maximal Jordan number $2m$, which form a single class under the action of $O^+(4m,2)$ \cite{ibw1}.   By conjugacy in $O^+(4m,2)$, $u$ lies in such a group, $P$.  
\eop

\begin{lem}\labtt{singer2}  
Suppose that $d\ge 2$ and that 
$t \in \gd d$ has defect $\frac d 2$ or $\frac {d-1} 2$.  Let $R:=\rd d$.  
Suppose that $[R,t]$ is elementary abelian.   Then 
$t$ is in a dihedral group as in \refpp{ineq2}.
\end{lem}
\pf  Let bars indicate images in $\gd d/\rd d$.   
Lemma \refpp{singer1} implies that $\bar t$ 
is in an appropriate Singer normalizer, $E$.
Let $u$ be a conjugate of $t$ in $G$ so that 
$\bar t \bar u$ generates 
$O_{2'}(E)$.  There exists $c \in \la tu \ra$ which generates a cyclic group of odd order which maps isomorphically onto $O_{2'}(E)$.  Then $\la t, c \ra$ satisfies the conclusion.  
\eop

Now we prove the main result \refpp{jno=elm}.

\begin{thm}\labtt{jno=elm}   Let $d\ge 2$ and let 
$t$ be an upper involution in $\brw d$ of defect $k\ge 1$.  
Then  $JNo(t)=\elm d k$ if $t$ is split, 
and is $2^{d-1}$ if $t$ is nonsplit.  
\end{thm}
\pf 
We have $d\ge 2$.  
Suppose that $[R,t]$ is not elementary abelian.  
There exists a lower involution $w$ so that $[w,t]$ has order 4.  Then on the 
lower dihedral group $D:=\la w, [w,t] \ra$,  $t$ induces an outer automorphism. Now use \refpp{nonsplitfull}.

We may assume that $t$ is split.  So, $[R,t]$ is elementary abelian.  If the involution $t$ centralizes a lower dihedral group, the 2/4 generation property \refpp{2/4}  and induction 
\refpp{eq1}  implies the result.  Note that the initial cases for induction are discussed in \cite{ibw1}.  

Assume that the involution $t$ does not centralize a lower dihedral group.  
Then  $\rd d/\rd d'$ is a free $\FF_2[\la t \ra ]$-module, $d$ is even and $d=2k$.  
We apply \refpp{singer2}, \refpp{singer1} with $r=m=k$, then \refpp{ineq3} and \refpp{ineq1}.  
\eop

\subsection{  Applications to discriminant groups  } 

Knowing $JNo$ is quite useful.  One can get sharp statements about the discriminant group, which might be hard to calculate directly from  a definition of the lattice, e.g. by a spanning set.    

\begin{lem}\labtt{commnegeigen}
Let the involution $u$ act on the additive abelian group $A$.  Then
$2A^- \le [A,u] \le A^-$.  
\end{lem}
\pf
Clearly, $u$ negates all $a(u-1)$, so $[A,u] \le A^-$.  
Also, if $a \in A^-$, $2a=a-(-a)=a(1-u)\in [A,u]$. 
\eop

\begin{coro}\labtt{minus=comm} Suppose that $L\cong \bw d$ and $t\in 
\gd d$ satisfies  $tr(t)>0$.  Then 
$L^-(t)=[L,t]$.  
\end{coro}
\pf  Since $t$ is an involution, $L^-(t)\ge [L,t]$ \refpp{commnegeigen}.  
Since $JNo(t)=rank(L^-(t))$ \refpp{jno=elm}, 
the image in $L/2L$ of $[L,t]$ has dimension equal to the rank of $[L,t]$.  
Therefore, $\lmt + 2L=[L,t]+2L$.   Since $[L,t]\le \lmt \le [L,t]+2L$, the Dedekind law implies that 
$\lmt \le [L,t]+(\lmt \cap 2L)$.  Since $\lmt$ is a direct summand of $L$, 
$\lmt \cap 2L=2\lmt$.  The latter is contained in $[L,t]$, by \refpp{commnegeigen}.  We conclude that $\lmt = [L,t]$.   
\eop


\begin{coro}\labtt{discriminantrank} 
Let $d\ge 2$.  
Let $t$ be a split involution of defect $k\ge 1$, and $\vep =\pm$.  Suppose $tr(t)>0$.  

(i)  The image of $L$ in the discriminant group of $L^\vep (t)$ 
is 2-elementary abelian of rank $\elm d k$. 

(ii) $\lmt \le 2P^-(t)$.  

(iii) 
If $d$ is odd, 
$\dg {L^-(t)}\cong \dg {L^+(t)}$ is 2-elementary abelian of rank $\elm d k$.  
In particular, ${L^- (t)} =2  \dual {L^- (t)} = 2P^-(L)$. 
\end{coro}
\pf 
(i)  The kernel of the natural map 
$\pi_\vep : L\rightarrow \dg {L^\vep(t)}$ is $L^+(t)\perp L^-(t)$.  The cokernel is elementary abelian of rank $JNo(t)$.  

(ii) Use (i) and rank considerations.  

(iii) Since $d$ is odd, unimodularity of $L$ implies that  each $\pi_\vep$ is onto.   
\eop

\section{Midwest cousins}

We introduce the first midwest operator here.  

\begin{de}
\labttr{mc}
The {\it midwest cousin (MC) lattices} are defined as follows.  
Let $L$ be an integral lattice.  
Let $t, f\in O(L)$ so that  $t, f$ commute, $t$ is an involution and $f$ is a fourvolution.   
Let $\vep = \pm$ 
and let 
$P^\vep$ be the orthogonal projection to $V^\vep (t)$.  
Set $MC(L,t,f, \vep ):=\lvept  + P^\vep (L)(f-1)=\lvept + P^\vep (L(f-1))$ (see \refpp{twist}\refpp{commdensapp}\refpp{commdensapp2} about alternate notation $L[p]$).  
\end{de}

\begin{lem}\labtt{mc2} Let $\vep = \pm$.   

(i) The midwest cousin $MC(L,t,f, \vep )$ 
is an integral lattice.

(ii) If $L^{\vep}(t)$ is doubly even, i.e., all norms are multiples of 4, then  
$MC(L,t,f, \vep)$ is an even lattice.  
\end{lem}
\pf
(i) 
We verify that $(x,y)\in \ZZ$, for $x, y \in MC(L,t,f, \vep )$.    
If $x$ or $y$ is in $\lvept\le L$, this is clear.  
Now suppose that $x=x'(f-1), y=y'(f-1)$ for $x', y' \in P^\vep (L)$.  
Then 
$(x,y)=(x'(f-1),y'(f-1))=2(x',y')=(x',2y')\in 
(P^{\vep}(L),\lvept )\le (L,\lvept )\le (L,L) \le \ZZ$.  

(ii)  We take $x\in L, y:=P^{\vep}(x)$.  
Then $2y=P^{\vep}(2x)\in P^{\vep}(Tel(L,t))=L^{\vep}(t)$ so that $2y \in L^{\vep}(t)$.  
We have $(2y,2y)\in 4\ZZ$ since by hypothesis, $L^{\vep}(t)$ is doubly even.  Therefore, $(y,y)\in \ZZ$ and so $y(f-1)$ has even norm.  
Since $L^{\vep}(t)$ is even, and $(P^{\vep}(L),L^{\vep}(t))\le \ZZ$, it follows that 
$MC(L,t,f, \vep )$ is even.  
\eop

\begin{de}
\labttr{mcbw} 
The {\it midwest first cousins of the Barnes-Wall lattices} are defined as follows.  
They are the MC lattices  with input lattice $\bw d$ and 
a pair $t, f$ as in \refpp{mc} 
where  $t$ is positive trace defect $k$ involution and 
$f \in C_R(t)$ is a lower fourvolution \refpp{ibw1review}.  
When $k<\frac d2$, such pairs are unique up to conjugacy in $\brw d$.  In this case, we use the briefer notation $MC_1(d,k,\vep)$ for $MC(\bw d,t,f, \vep )$.  
When $k=\frac d2$, there are several conjugacy classes of pairs $(t,f)$.  One would need additional notation to distinguish these classes \cite{bwy}.  
\end{de} 

\begin{rem} 
\labttr{mcbw2} 
Let $L:=\bw d$.  
Suppose that we have two  pairs $(t, f)$ and $(t, f')$, 
where both $f, f'$ are lower fourvolutions which commute with $t$, 
then the resulting first cousin lattices are the same.  
The reasons are that $L(f-1)^p=L(f'-1)^p$, for all $p$ (because any lower fourvolution is commutator dense for the action of $R$ on $L$ \cite{bwy}) and the projection maps $P^\vep$ commute with $f$ and $f'$.  In certain commutator calculations, it may be convenient to replace $f-1$ by some $\pm f' \pm 1$.  
\end{rem}

\subsection{Integrality properties of the first cousin lattices}  

We now specialize to the case of Barnes-Wall lattices.  

\begin{prop}\labtt{propmcbw}  Let $d\ge 2$, $L:=\bw d$.  
We assume that the involution $t$  has defect $k\ge 1$ and that its trace is positive.   Then 

(i) $rank(MC_1(d, k, \pm) )=2^{d-1} \pm  2^{d-k-1}$.

(ii)  Let $\vep = \pm$.  
If $d$ is odd and $d\ge 3$, $MC_1(d,k,\vep)$ is unimodular.

(iii) For $\vep = \pm$,   
$k\le \frac d2 -1$, then $P^{\vep}(t)(f-1)$ is even integral and 
$\lvept$ is doubly even (and so $MC_1(d,k,\vep)$ is even).  

(iv) $\mu (MC_1(d,k,-))=\half \mu (\bw d )$.

(v) $\mu (MC_1(d,k,+)) \le  2^{\lfloor \frac d2 \rfloor}$.

(vi) If $d=2k$ or $d=2k+1$, $MC_1(d,k,\vep )$ is an odd integral lattice.  \end{prop}
\pf 
For (i), see \refpp{traces}.

For (ii), we have that $\half \lmt = P(L)$, which is $\dual{\lmt}$ since $L$ is unimodular \refpp{discriminantrank}(iii).  Consequently,  $\dg {\lmt}\cong 2^{rank(\lmt)}=2^{2^{d-1}-2^{d-k-1}}$.  The lattice $MC_1(d,k,-)$ is between $\lmt$ and its dual and corresponds to the image of $f-1$, where $f$ is a lower fourvolution in $C_R(t)$.  
In fact, $MC_1(d,k,-)=P^-(L)(f-1)$.  
Since $(f-1)^2=-2f$ and $|\half \lmt : MC_1(d,k,-)|=|MC_1(d,k,-):\lmt|$, unimodularity follows.  

The argument for $\vep = +$ is similar since 
$\dg {\lpt}\cong \dg {\lmt}$ as modules for $f-1$.

(iii) By \refpp{mc2}, $P^{\vep}(L)(f-1)$ is integral.  We show that it is even under our restrictions on $k$.  

Since $k<\frac d2$, there exists a lower dihedral group $D\le C_R(t)$ so that 
$D\cap [R,t]=Z(R)$.  If $u, v$ form a generating set of involutions, $L=L^+(u)+L^+(v)$ by 2/4-generation \refpp{2/4}.  The action of $t$ on each summand has nonzero trace and defect $k$.  

Suppose that $d$ is even.   Then $d-1$ is odd and each summand is $t$-invariant and is isometric to $\sqrt 2 \bw {d-1}$.  By a previous paragraph, the norms of vectors in $P^{\vep}(L^+(u))$ and $P^{\vep}(L^+(v))$ are integral.  Therefore 
the norms of vectors in $P^{\vep}(L^+(u))(f-1)$ and $P^{\vep}(L^+(v))(f-1)$ are even integral.  
This suffices to prove (iii) since we have a spanning set of even vectors in an integral lattice.

For (iv), note that $L^-(t)$ contains a minimal vector of $L$  and that $MC_1(d,k,-)$ is the $-1$ twist \refpp{twist} of $L^-(t)$.  

(v) This is obvious since $\lpt$ contains a minimal vector of $L$.   

(vi) 
Integrality was proved in \refpp{mc2}(i).  

If $d=2k$, the vector $v:=2^{-k}v_Z$ is in $P^{\vep}(L)$.  Its norm is 
$2^{-2k}2^k(2^{d-1}+2^{d-k-1})=2^{d-k-1}+\vep \half$.  The vector 
$v(f-1)$ is in $MC_1(d,k,+)$ and has odd integer norm.

If $d=2k+1$, let $H$ be an affine hyperplane which is transverse to $core(Z)$, which is 1-dimensional.  
The vector $v:=2^{-k}v_{H\cap Z}$ is in $P^+(L)$ 
and has norm 
$2^{-2k}2^k(2^{d-2}+2^{d-k-2})=2^{d-k-2}\vep \half$.
The vector 
$v(f-1)$ is in $MC_1(d,k,-)$ and has odd integer norm.  
To prove the result for $\vep=-$, replace $Z$ by $Z+\Omega $ in the above reasoning.

Suppose that $d=2k$ is even.  Then $2^{-k}v_{\Omega}\in L$ and
$2^{-k}v_{Z}\in P^+(L)$.  Its norm is 
$2^{-2k}2^k(2^{d-1}+\vep 2^{d-k-1})=2^{k-1}+\half$.  
The vector 
$v(f-1)$ is in $MC_1(d,k,+)$ and has odd integer norm.  A similar argument works for $\vep = -$.  
\eop 

\begin{rem}\labttr{mc52}
The unimodular 
integral lattices $MC_1(5,2,\pm)$ are  not even since their ranks are 20 and  12, which are not multiples of 8.  Another proof is \refpp{propmcbw}.  
\end{rem}

\subsection{Minimum norm for $MC_1(d,k,+)$}

In this section,  we determine that the minimum norm for $MC(d,k,+)$ is $2^{\frac {d-1}2-1}$ \refpp{minvecdescr}, the same as for $MC_1(d,k,-)$\refpp{propmcbw}.  Later, we discuss the forms for low norm vectors in the first few layers \refpp{layers} and study orthogonal decomposability.

\begin{nota}\labttr{notatop} We let $t$ be an involution of defect $k$ and positive trace.  We take $t$ to have the form $\vep_{Z}$, 
where $Z$ 
has weight $2^{d-1}+2^{d-k-1}$.    
As before, abbreviate $P^\vep$ for the projection to $L^\vep (t)$.  
Let $c\in core(Z), c\ne 1$ \refpp{core} and let $H$ be a hyperplane of $\Omega$ which is transverse to $\{0,c\}$ (so is moved by translation by $c$).  
We take $\tau:=\tau_c$, $f:=\vep_H\tau$ and define $\xi :=f-1$, so that $L[k]=L \xi^k$, for all $k$.  
\end{nota}

\begin{nota}\labttr{delta}  $\delta:=\frac {d-1}2$.  
\end{nota}

\begin{thm} \labtt{minvecdescr}  
We suppose that $d-2k\ge 3$.  

(i) 
$\mu (MC_1(d,k,\vep )) = 2^{\delta -1}$. 

(ii) A vector  $v\in MC_1(d,k,\vep )$ is minimal if and only if $v\xi$ is minimal in $\lvept$ (equivalently, if the support of $v\xi$ is contained in $Z$ 
and $v\xi$ is a minimal vector of $\bw d$).  

(iii) The minimal vectors of $MC_1(d,k,\vep )$ are in $MC_1(d,k,\vep ) \setminus \lvept$.  
\end{thm}
\pf
(i) Let $v\in \mv {MC_1(d,k,\vep )}$.  Since $v\xi \in \lpt$, $(v,v)\ge 2^{\delta -1}$.  It suffices to prove that there exists a vector in $MC_1(d,k,\vep )$ of such a norm.

We let $p\ge 1$ and let $A$ be an affine subspace of dimension $2p$  in $\Omega$ which is 
a translation of a subspace of $core(Z)$ (this is possible since 
$d-2k\ge 3$).  We also choose $A$ to be transverse to $H$ (this is possible since $2p<d-2k$) and to be contained in $Z$.  Therefore, $A\cap H$ is a $(2p-1)$-dimensional space.  
The vector $2^{-p}v_{A\cap H}$ is in $MC_1(d,k,\vep )$ and has norm $2^{\delta -1}$.  

(ii) Since $\xi$ takes $MC_1(d,k,\vep )$ into $\lpt$ and doubles norms, this follows from (i).  

(iii) 
This follows from (ii) since the minimum norm in $L$ is $2^{\delta}$.  
\eop

\begin{coro}\labtt{minvecdescr1} 
A minimal vector of $MC_1(d,k,\vep )$ has the form 
$2^{-m}v_A\vep_S$, where $A$ is an affine $(2m-1)$-space, $A\subseteq Z$ and $S\in RM(2,d)$.  
\end{coro}
\pf
Use  \refpp{minvecdescr}(ii), \refpp{minvecbw}, \refpp{minvecbwd[1]}.  
\eop

\begin{rem}\labttr{minvecdescr2} 
The description \refpp{minvecdescr1}  of minimal vectors in $MC_1(d,k,\vep )$ is similar to \refpp{minvecbw} for $\bw d$, but is not as definitive.   
\end{rem}

\section{Lattices with binary bases}  

To prove our main results about short vectors in the  lattices $MC_1(d,k,\vep)$, we begin with a general theory for  lattices with a binary basis.   Later, we shall specialize to the Barnes-Wall lattices.

\begin{de}\labttr{layers} 
Let $L$ be an integral lattice and  $M$ in another lattice in $\QQ \otimes L$ so that 
$L\le \ZZ [\half ]\otimes M$.  Let $q\ge 0$ be an integer.  
Define $L(q):=2^{-q}M \cap L$.  Call this the {\it $M$-level $q$ sublattice of $L$}. 
The {\it level} of $0\ne x\in L$ with respect to $M$ is $min\{k \ge 0 \mid x\in L(k)\}$.  
The {\it $q$-th layer} of $L$ is $L(q)/L(q-1)$.  If $S$ is a subset of $\QQ \otimes L$ which is $\QQ$-linearly independent and such that its $\ZZ [\half]$-span contains $L$, we call $S$ a {\it binary basis} and define 
 {\it level} of $x\in L$ with respect to $S$ to be 
the level of $x\in L$ with respect to $span_{\ZZ}(S)$.  We do not assume that $S$ is an orthogonal set.  
\end{de}

 \begin{nota}\labttr{2adic} If $n\in \ZZ [\half]$ is nonnegative, its {\it 2-adic expansion} is an expression 
 $n=\sum_{i=p}^q a_i2^i$, where the $a_i$ come from $\{0,1\}$.  When $n\in \ZZ [\half]$ is negative, its 2-adic expansion is $\sum_{i=p}^q -a_i2^i$, where $-n=\sum_{i=p}^q a_i2^i$ is the 2-adic expansion of the nonnegative rational $-n$.  
 The {\it level of $n$} is $-\infty$ if $n=0$  and is otherwise $-min\{i \mid a_i\ne 0\}$.  
 \end{nota}

 \begin{nota}\labttr{top} Let $L$ be a lattice of rank $n$ with $S$, a linearly independent subset 
$v_1, \cdots, v_n$.  Then $x \in L$ has a unique expression 
$x=\sum_i c_iv_i$, for rational numbers $c_i$.  We assume that 
$S$ is a  binary basis for $L$ \refpp{layers}.   
  Then the $c_i$ are in $\ZZ [\half ]$. 

We define the {\it 2-adic expansion of $x$} to be $\sum_i 2^i (\sum_j a_{i,j} v_j)$ where the $a_{i,j}$ are the 2-adic coefficients of $c_j$.  For $x \in L$, define $level(x)$ to be the least integer $m$ so that 
the coefficients of 
$\sum_i 2^m c_iv_i$ 
are integers.  We define $level(0):=-\infty$. 

For $x\ne 0$, we define $top(x)=top_S(x)$ to be the subsum  $\sum_j a_{m,j} v_j$ 
of the 2-adic expansion of $x$ (it is the part of the 2-adic expansion of $x$ which represents the largest denominators, $2^m$).   Note that the definition of $top(x)$ depends on the binary basis, not on the sublattice it spans.  
\end{nota}

 \begin{rem}\labttr{topnotin} (i) The top of a vector may not be in the lattice.  Consider the lattice $L$  in $\QQ^2$ which is spanned over $\ZZ$ by $(1,0), (0,1), (\half, \fourth)$.  For $S$, take $\{(1,0), (0,1)\}$.  We claim that  $top((\half, \fourth))=(0,\fourth)$ is not in $L$.  If $(0,\fourth)=a(1,0)+b(0,1)+c(\half,\fourth)$, we may assume that $c\in \{0,1,2,3\}$.  Clearly, $c$ is $1(mod \ 4)$, so $c=1$.  Then the right side has first coordinate a noninteger, contradiction.  
 
 (ii) Tops do lie in $\bw d$ for  
 vectors of level at most 1 with respect to the 
 the standard basis in a lower frame.  For higher level, top closure may fail.  For example, take $d\ge 8$ and consider a pair of 4-spaces which meet in a point.  
 \end{rem}

\section{Calculations in $MC_1(d,k,\vep )$}

\begin{coro}\labtt{topinmc}
Suppose that $0\ne x\in MC_1(d,k,\vep )$ has level $m$.  
Then 
$top(x)=2^{-m}v_B$, where $B\in RM(d-2m+1,d)$.  
Furthermore, given $\tau = 
\tau_c$ in $0\ne c \in core(Z)$, there is a decomposition $B=S+T$, where 

(i) $S\in RM(d-2m,d),  T\in RM(d-2m+1, d)$; 

(ii) $S\subseteq Z, T\subseteq Z$;  and 

(iii) $T$ is $\tau$-invariant or $T$ has form $A\cap H$ 
where $A\in RM(d-2m+2,d)$, $A\subseteq Z$, $A$ is $\tau$-invariant and $H$ is a hyperplane transverse to $\tau$ (i.e., transverse to $\{0,c\}$ in $\Omega$).  
\end{coro} 
\pf 
Since $MC_1(d,k,\vep )=\lvept + P^{\vep}(t) [1]$, this 
follows from 
the corresponding forms for $top(x)$, $x\in \lvept$ and $x\in P^{\vep}(t) [1]$ and the action of $f-1$.  
\eop

\subsection{Equations with codewords and commutation}

We collect a few results about expressions of the form 
$B=S+T\in RM(d-2m,d)$ as in \refpp{topinmc}.

\begin{lem}\labtt{zdec2} Suppose that $B\in RM(i,d)$, 
$B=S+T\in RM(d-2m,d)$ as in \refpp{topinmc}.  
Let $r$ be a real number so that $|B|\le 2^r$.  
If $d> r+i$, then $B$ is $\tau$-invariant.
\end{lem} 
\pf   
We may assume that $i\ge 1$.  
We have $B(\tau -1)\in RM(i-1,d)$, which has minimum weight $2^{d-(i-1)}$.  
Since $|B(\tau -1)|\le 2^{r+1}$, if $B(\tau -1)\ne 0$, then $d-i+1\le r+1$, or $d\le r+i$, contrary to hypothesis.  Therefore $B(\tau -1)=0$, i.e., $B$ is $\tau$-invariant. 
\eop

\begin{coro}\labtt{r=1}  Assume the hypotheses of   \refpp{zdec2}.   If $0\ne |B|\le 2$ and $i=d-2$, then $B$ is  $\tau$-invariant.
\end{coro}
\pf Take $r=1$ in \refpp{zdec2}.  
\eop

\begin{lem}\labtt{bfixed} Suppose $\tau =  \tau_c$, for $c\in core(Z)$ and $c\ne 0$.  
 Suppose $B\in RM(d-2m+1,d)$ is fixed by $\tau$. 
Then $|B|\ge 2^{2m-1}$.  
\end{lem}
\pf  Let bars denote images in the quotient code $\Omega / \Gamma$ \refpp{quotientcode}, where $\Gamma = \{0,c\}$.   
Then  $\bar B$ is a nontrivial element of $RM(d-2m+1,d-1)=RM((d-1)-(2m-2),d-1)$, so has weight at least $2^{2m-2}$.  This implies $|B|\ge 2^{2m-1}$.  
\eop

\section{$MC_1(d,k,\vep)$ short vectors, level at most 2}

By \refpp{minvecdescr1}, a minimal vector of $MC_1(d,k,\vep)$ is a vector of the form $2^{-m}v_B\vep_C$, for some $m\ge 0$, some $B\in RM(d-2m+1,d)$ and some $C\subseteq \Omega$.  We can say more about short vectors in the first two levels.  

Recall the concept of level \refpp{layers}.  
Vectors of level 0 are in $\bw d$, so their norms are 0 or are at least $2^\delta$.  The set of level 0 norm $2^\d$ vectors is just $\{\pm v_i  \mid  i \in \ZZ\}$, the standard lower frame.  

\subsection{Short vectors at level 1}

We display a set of norm $2^{\d- 1}$ vectors, which turn out to be the only 
level 1 vectors in $MC_1(d,k,\vep )$ of norm less than $2^\delta$.

\begin{lem} \labtt{level1candidates}  Suppose that the level of $0\ne x \in MC_1(d,k, \vep)$ is 1.  So, $top(x)=\half v_B$.  Then:

(i)  $|B|$ is even.  

(ii) 
If $(x,x)<2^\delta$, 
then $B$ is a 2-set and  $B$ is stabilized by some $\tau_c\ne 1$.       
\end{lem} 
\pf 
(i) Trivial 
since $B\in RM(d-2m+1,d)$ and $m=1$.  

(ii) Use \refpp{r=1}.    \eop

\begin{lem} \labtt{level1list}  The set of level 1 vectors of $MC_1(d,k,\vep )$ of norm less than $2^\delta$ consists of all 
$\pm \half v_i \pm \half v_{i+c}$, for $c\ne 1, c\in core(Z)$ and $i \in \Omega$.  These have norm $2^{\d-1}$.   
\end{lem}
\pf  We get a list of candidates from \refpp{level1candidates}(ii).  
We need to see that all the vectors of indicated form are actually in $MC_1(d,k,\vep )$. 
By \refpp{rmduality}, there exists $E\in RM(d-2,d)$ so that $F:=E\cap Z$ is an odd set.  Therefore $F(\tau -1)$ has cardinality $2(mod\ 4)$.  By \refpp{auf}(ii), there exists $S\in RM(d-2,d)$ so that $B=S+F(\tau -1)$ is a 2-set, and such a 2-set is $\tau$-invariant \refpp{r=1} and so is one of the indicated $\{i, i+c\}$.  

\begin{rem}\labttr{graph}  We recall an elementary result about positive definite integral lattices \cite{kn}.  Let $J$ be such a lattice.  Call $x\in J, x\ne 0$ {\it decomposable} if there exist nonzero $y, z\in J$ so that $x=y+z$.  If $X$ is the set of indecomposable vectors, we define a graph structure by connecting two members of $X$ with an edge if they are not orthogonal.  We therefore get $X$ as the disjoint union of connected components $X_i$.  If $J_i$ is the sublattice spanned by $X_i$, then $X$ is their orthogonal direct sum.   If $Y$ is any orthogonal direct summand of $J$, $Y$ is a sum of a subset of the $J_i$.  
\end{rem}

\begin{coro}\labtt{level1span}
The vectors of \refpp{level1list} span a sublattice which is an orthogonal direct sum of 
scaled $D_{2^{d-2k}}$ root lattices.  This sublattice has finite index in $MC_1(d,k,\vep )$.  
\end{coro} 
\pf
Consider the natural graph on this set of vectors where edges between distinct vectors are based on nonorthogonality.  The connected components span lattices of type $D$ \refpp{level1list}.  
\eop

\subsection{Short vectors at level 2}

For the moment, $d\ge 5$ is odd and arbitrary.  Recall that 
top closure may fail in $\bw d$ above level 1 \refpp{topnotin}.

\begin{prop}\labtt{lev2def1d=5}  Suppose that $d\ge 5$ and $d-2k\ge 3$.   
If the norm of the level 2 vector $x\in MC_1(d,k,\vep )$ is $2^{\d -1}$,  then there exists 
$C\in \pow {\Omega}$ and $B$ is affine 3-space so that 
$x=\fourth v_B\vep_C$.   
\end{prop} 
\pf   
Since $B\in RM(d-2m+1,d)$, we use \refpp{minweight}.
\eop

\begin{rem}\labttr{lev2def1d=5} We do not assert that vectors as in 
\refpp{lev2def1d=5} exist. 
\end{rem}

\section{Decomposability  and indecomposability}

We prove that the first cousins are orthogonally decomposable for $k=1$ and indecomposable for $k\ge 2$.  As in \refpp{notatop}, $t$ has positive trace.  

\begin{prop}\labtt{mck=1eps=-}  Let $k=1$.  
The lattice $MC_1(d,1,-)$ is isometric to $\bw {d-2}$.  
\end{prop}
\pf
By ancestral theory \cite{bwy}, $\lmt \cong \bw {d-2}[1]$.  By \refpp{discriminantrank}(iii), 
$MC_1(d,1,-)\cong 2^{-\half}\lmt \cong \bw {d-2}$.
\eop 

\begin{prop}\labtt{mck=1}  Let $k=1$.  
The lattice $MC_1(d,1,+)$ is isometric to $\bw {d-2} \perp \bw {d-2} \perp \bw {d-2}$.  
\end{prop}
\pf
By hypothesis, $k=1$.  Thus, $Z$ is the complement in $\Omega$ of a codimension 2 affine space.  
There are three affine hyperplanes contained in $Z$.   Call them $Z_1, Z_2, Z_3$ and let $Z_{ij}$ denote the intersection of $Z_i$ and $Z_j$. 

The proof is a consequence of the theory of \cite{bwy}.  
For a subset $T$ of $\Omega$, we let 
$L(T)$ be the set of vectors in $L$ whose support is contained in $T$.  Then $L(Z_i)$ is a scaled $\bw {d-1}$.  
The sublattice $L(Z)$ is coelementary abelian of 
index $2^{2^{d-2}}$ in the orthogonal direct sum 
$\half L(Z_{12}) \perp \half L(Z_{23}) \perp \half L(Z_{31})$.   
Furthermore, a set of coset representatives for $L(Z_i)\perp L(\Omega + Z_i)$ in $L$ is just 
the set $S$ of all $x+xu$, where $u$ is a fixed involution interchanging $L(Z_i)$ and $L(\Omega + Z_i)$ and where $x\in L(Z_i)[-1]$. 
(The relevant lower fourvolution $f$ should be chosen to have an expression $f=\prod f_i$, where $f_i$ is a lower fourvolution on $Z_i$; see \cite{bwy,ibw1}).  
 
It follows that the set $P^{+}(S)(f-1)$ represents all the cosets of $L(Z)$ in 
$\half L(Z_{12}) \perp \half L(Z_{23}) \perp \half L(Z_{31})$.  (It may help to think that $\FF_2^3$ is spanned by $(1,1,1)$ and the space of vectors with coordinate sum 0.) 
\eop

\begin{lem}\labttr{indecfinindex} 
Suppose that $M$ is an integral lattice and $N$ a finite index sublattice.   Suppose that $N$ is spanned by vectors which are indecomposable in $M$ and that $N$ is orthogonally indecomposable.  Then $M$ is orthogonally indecomposable. 
\end{lem} 
\pf  
 The hypotheses on $M$ and $N$ imply that $N$ meets every indecomposable summand of $M$ nontrivially.  See \refpp{graph}.  
\eop 

\begin{lem}\labtt{highestlevelminvec}
Recall that $H$ is a hyperplane which is transverse to $core(Z)$.  Set $v:=2^{-\d}v_H$, a minimal vector in $\bw d$.  
Then $P^{\vep}(v)$ has norm $2^{\d} \frac {|Z|}{2^d}=r2^{\d}$, for some $r\in [\frac 14, \frac 34]$.  Also, 
$P^{\vep}(v)(f-1)$ has norm $r2^{\d+1}=s2^{\d -1}$, for some $s\in [1,3]$.  Therefore, if we write 
$P^{\vep}(v)(f-1)=w_1+\dots + w_n$ as an orthogonal sum of indecomposable nonzero vectors, $n\le 3$.  
\end{lem}  
\pf 
Use
the formula for $|Z|$ \refpp{notatop}, 
\refpp{minvecdescr} and the fact that 
$P^{\vep}(v)(f-1) \in MC_1(d,k,\vep )$.  
\eop

\begin{prop}\labtt{lev2def1dge7}  
Suppose that $d\ge 7$ is odd and $k\ge 2$.  

(i) 
The minimal vectors of the level 1 sublattice are indecomposable in $MC_1(d,k,+ )$.  The sublattice of $MC_1(d,k,+ )$ which they span is an orthogonal direct sum of scaled type $D_{2^{d-2k}}$ lattices. 

(ii) 
When $d\ge 7$ and $d-2k\ge 5$, the lattice 
spanned by the level 2 minimal vectors (which have norms $2^{\delta - 1}$) is orthogonally indecomposable and has finite index in $MC_1(d,k,+ )$.  
Therefore,  $MC_1(d,1,+)$ is orthogonally indecomposable.  
\end{prop}
\pf 
(i) 
The first statement is trivial since they are minimal vectors in $MC_1(d,k,+ )$.  The second statement follows from analysis as in the proof of \refpp{lev2def1d=5}. 

(ii) Let 
$L_1, \dots, L_r$ be the set of scaled type $D_{2^{d-2k}}$-lattices as described in (i).  Each is orthogonally indecomposable since $d-2k\ge 3$.  

Take a vector hyperplane $H$ and vector $v$ as in \refpp{highestlevelminvec}.  Then $v$ has nonzero inner product with vectors of each $L_i$ and so does $P^{+}(v)(f-1)$.  If we write $P^{+}(v)(f-1)=w_1+\dots + w_n$ as a sum of indecomposable vectors, we get $n\le 3$ by norm considerations.  For each $i$, there exists $j$ so that $L_i$ has nonzero inner products with $w_j$.   
The number of $L_i$ is $2^{d-1}+2^{d-k-1}$, which is at least 4, and the number of $w_j$ is at most 3.   
Therefore, there exists a pair of distinct indices $i, i'$ and an index $j$ so that
both $(L_i,w_j)$ and $(L_{i'},w_j)$ are nonzero.   Therefore in the graph of indecomposable vectors \refpp{graph}, the minimal vectors of $L_i$ and $L_{i'}$ are in the same component.  Now we quote double transitivity of $Sp(2k,2)$ on the set of $L_i$ \cite{ibw1} to deduce that all minimal vectors of $L_1\perp L_2\perp \dots \perp L_r$ are in the same component.  This proves that $MC_1(d,k,+ )$ is indecomposable. 
\eop

\section{More distant cousins}

We have considered variations of the formula for first cousins.  
Many interesting high dimensional lattices with moderately high minimum norms may be created in the midwest style.  Precise analysis of their properties would be challenging, however.  

One variation creates an even unimodular rank 24 overlattice of $\lpt$ for $L\cong \bw 4$ and $tr(t)=8$.  That overlattice has minimum norm 4, so is isometric to the Leech lattice.

Here is a sketch of the construction.   In $\lpt$, there is a sublattice $M=M_1\perp M_2\perp M_3$, where  $M_i\cong \sqrt 2 E_8$, for $i=1,2,3$.  Let $f$ be a lower fourvolution on $L$ which commutes with $t$ and fixes each $M_i$.  Then $\lpt (f-1)\le M$ and $P^+(L)(f-1)\le \lpt$.  
We need a lemma.  

\begin{lem}\labtt{gexists}
Suppose that we have two sublattices 
$M, N$ such that $E_8=M+N$ and $M\cong N\cong \sqrt 2 E_8$.
There exists $\g \in O(E_8)$ which interchanges $M$ and $N$.  
\end{lem} 
\pf
This follows from the analogous property of $O^+(2d,2)$ since $O(E_8)$ acts on $E_8$ mod 2 as $O^+(8,2)$.  
\eop

\smallskip

Continuing our construction, we let $\g$ be an isometry of $M$ which stabilizes each $M_i$ and satisfies $M_i(f-1) \cap M_i(f-1)\g = 2M_i$ and (consequently) that 
$M_i(f-1) + M_i(f-1)\g = M_i$ (see \refpp{gexists} and the ancestral theory \cite{bwy}).   Then 
$\lpt + P^+(L)(\g^{-1} f\g -1)^2$ is isometric to the Leech lattice.  There is similarity in spirit to \cite{lm,tits}.  

It is well-known that the Leech lattice contains sublattices isometric to $\bw 4$ (as fixed point sublattices of involutions) \cite{cs}, \cite{poe}.  The above result links the Leech lattice and $\bw 5$.

\section{Appendix: Some  background} 

Standard properties of Reed-Muller binary codes \cite{reed, muller} and the Barnes-Wall lattices \cite{bw, be, bwy}   will be used intensely.  For convenience, we review them here.  

\subsection{Review of Reed-Muller codes}

\begin{nota}\labttr{rmsetup} 
For integers $d\ge 1$ and $k\in \{0, 1, \cdots , d\}$, there is defined a Reed-Muller binary code $RM(k,d)$ of length $2^d$.   We use $\Omega=\Omega_d$, a copy of affine space $\FF_2^d$, as indices.   A binary vector may be interpreted as an $\FF_2$-valued function of its index set $\FF_2^d$, or as a subset of the index set (the support of the previous function).  Addition is the boolean sum.  The Reed-Muller code $RM(k,d)$ is spanned by the vectors which are the characteristic functions of affine subspaces of codimension at most $k$ (or, in the power set interpretation $\FF_2^{\Omega}$, as the actual affine subspaces).  For all $p\le -1$, $RM(p,d):=0$.  
\end{nota}  

We mention a few facts for use in this article.

\begin{prop}\labtt{rmduality} For $d\ge 1$ and  for $i=0,1,2, \dots , d-1$, $RM(i,d)^\perp = RM(d-1-i,d)$.   
\end{prop} 

\begin{lem}\labtt{minweight}  
In $RM(k,d)$, the minimum weight is $2^{d-k}$ and the codewords of minimum weight are the affine subspaces of codimension $k$;

\end{lem}
\pf 
This is well-known; see \cite{ms}, Theorem 3, p. 375 and Theorem 8, p. 380.  

\eop

\def\rmlev#1{\text{RM-level}(#1)} 
\def\bwlev#1{\text{BW-level}(#1)}

\begin{de}\labttr{levels}
For $A\in \pow {\Omega}$, we define the {\it BW-level} of $A$ to be $max\{m\ge 0 \mid  A\in RM(d-2m,d)\}$ and the {\it RM-level} of $A$ to be
$max\{i \mid A \in RM(d-i,d)\}$.  We abbreviate these 
terms by {\it BW-level($A$)} and  {\it RM-level($A$)}, respectively.  
We extend the concept of level to elements of $\bw d$ by using the notation \refpp{top} with respect to the basis $v_i$ of \refpp{rmsetup}.  
\end{de} 

\begin{rem}\labttr{levels2}
If $i=\rmlev A$, then the elements of 
$A+RM(d-i-1,d)$ have RM-level $i$.  
If $m=\bwlev A$, then the elements of
$A+RM(d-2m-2,d)$ have BW-level $m$.  
\end{rem}

\begin{prop}
\labtt{tau-1}
Suppose that $\tau $ is a translation in $AGL(d,2)$.  Then 

(i) $RM(j,d)(\tau -1)\le RM(j-1,d)$; 

(ii) $\pow {\Omega}$ is a free module for 
$\FF_2 [\FF_2^d]$.  
The image of $\tau -1$ is the set of all $\tau$-invariant codewords. 
Also,  $\pow {\Omega}$ is a free $\FF_2 [\la \tau \ra]$-module.  

(iii) If $x \in Ker(\tau - 1)=Im(\tau -1)$ and $x\in RM(d-k,d)$, there exists  $y \in RM(d-k+1,d)$ so that $x=y(\tau -1)$.

(iv) If we identify the group algebra $\FF_2 [\FF_2^d]$ with $\pow {\Omega}$, the powers of the augmentation ideal of $\FF_2 [\FF_2^d]$ are the codes $RM(j,d)$.  
\end{prop} 
\pf  
(i) The first part is obvious since $RM(j,d)$ is spanned by affine subspaces $S$ of codimension $j$, and $S+S\tau$ is either empty or is a 
$(j+1)$-dimensional affine subspace.  

(ii)  Since $\pow {\Omega}$ is a free module for 
$\FF_2 [\FF_2^d]$, it is a free module for the subalgebra 
$\FF_2 [\la \tau \ra]$.  The statements follow.

(iii)  Since $\pow {\Omega}$ is a free module for  $\FF_2 [\la \tau \ra]$ (by (ii)), $Ker(\tau -1)=Im(\tau -1)$.    Assume that $c$ is a $\tau$-invariant codeword in $RM(k,d)$.   Since $\tau$ is an involution, $c$ is an even set, whence $k \le d-1$.    Let $h$ be an affine hyperplane which is transverse to every $\tau$-invariant 1-space.  Then $c\cap h\in RM(k+1,d)$ and $c=(c\cap h)(\tau -1)$.  

(iv) This follows from (ii) and (iii).  
\eop

\begin{lem}\labtt{auf}
Let $X$ be a subset of $\Omega$.   Then 

(i) if $|X|$ is even, $X(\tau -1)$ is in $RM(d-2,d)$; and

(ii) if $|X|$ is odd,  there is $Q$, a 1-space invariant under $\tau$,  such that $X(\tau -1)$ is in $Q+RM(d-2,d)$.

(iii) In (ii), if $Q, Q'$ are 1-spaces such that $X(\tau -1)$ is in $Q+RM(d-2,d)=Q'+RM(d-2,d)$, then $Q'$ is a translate of $Q$ and \ both are $\tau$-invariant. 
\end{lem}
\pf  
To prove (i), use \refpp{tau-1}(i).  Next,  (ii)  follow easily from the case $|X|=1$.  
For (iii), we may assume $X$ is a 1-set.  First notice that since $Q+Q'\in RM(d-2,d)$, whose minimal weight codewords are affine 2-spaces, $Q'$ is a translate of $Q$.  One is $\tau$-invariant if and only if the other one is.  On the other hand, there exists some 1-space $Q''$ which is $\tau$-invariant and which satisfies $X(\tau -1)\in Q''+RM(d-2,d)$ (just take $Q''=\{x, x\tau\}$, for any $x\in X$, and use (i),(ii)).  Therefore, both $Q$ and $Q'$ are $\tau$-invariant.  
\eop

\begin{de}\labttr{quotientcode}
Suppose that $\Gamma$ is a subspace of $\Omega$.  Let $\pow{\Omega,\Gamma}$ be the members of $\pow{\Omega}$ which are unions of cosets of $\Gamma$.  Then members of ${\mathcal P}(\Omega,\Gamma)$ may be interpreted as subsets of the quotient vector space $\Omega/\Gamma$ and so we have an isomorphism 
${\mathcal P}(\Omega,\Gamma)\rightarrow {\mathcal P}(\Omega/ \Gamma)$.  This may be interpreted as an isomorphism of a subspace of binary vectors of length $|\Omega|$ with the full space of binary vectors of length $|\Omega/\Gamma|$.
\end{de}

\begin{de}\labttr{termsrm2d} 
Given a codeword $c\in RM(2,d)$, there is at most one integer 
 $k\in \{1,2 \dots ,\frac d2 \}$  such that  the coset 
$c+RM(1,d)$ contains a codeword of weight  $2^{d-1}-2^{d-k-1}$.  
If there is such a $k$, we say $c$ has {\it defect $k$}.  If there is no such $k$, we say that $c$ has defect 0.  
We say that $c$ is {\it short} if it has cardinality less than $2^{d-1}$, {\it long} if it has cardinality greater than $2^{d-1}$ and otherwise we say $c$ is a {\it midset} or a {\it midword}.
\cite{ibw1} 
\end{de} 

\begin{de}\labttr{core} A sum $S_1+\cdots +S_k$ of $k>0$ affine codimension 2 subspaces whose intersection is nonempty, is called a {\it cubi sum} if its cardinality is $2^{d-1}-2^{d-k-1}$.  A short defect $k$ codeword $c$ may be written as a cubi sum.  We define the {\it core} of a cubi sum to be the intersection of the $k$ summands.   It depends only on $c$ and not on the particular cubi sum for $c$. 
\end{de}

\subsection{Review of PO$2^d$-theory and Barnes-Wall lattices} 

The Reed-Muller codes can be used to construct Barnes-Wall lattces \cite{bw}, \cite{be}.  Alternatively, they may be deduced from existence of Barnes-Wall lattices \cite{bwy}. 

\begin{nota}\labttr{bw0}
The {\it Barnes-Wall lattice} $\bw d$ in rank $2^d$, $d\ge 2$,  is an even 
lattice whose isometry group contains $\gd d \cong \ratholoex d$.  This is the full isometry group when $d\ne 3$.  These lattices are scaled so as to make $\bw d$ unimodular when $d$ is odd and to make the discriminant group elementary abelian of rank $2^{d-1}$  when $d$ is even.  Finally, define $\rd d:=O_2(\gd d)\cong 2^{1+2d}_+$.  
\end{nota}

 \begin{de} \labttr{standardgenset} 
 For $\bw d$, there is a standard generating sets (as abelian groups).    We start with the a set $\{v_i \mid i \in \Omega \}$ of vectors in $\bw d$.   As in \refpp{rmsetup}, $\Omega = \FF_2^d$.  We often use the maps $\vep_S$, which take $v_i$ to $-v_i$ if $i\in S$ and to $v_i$ if $i \not \in S$.  This map is in $\gd d$ if and only if $S\in RM(2,d)$ and is in $\rd d$ if and only if $S\in RM(1,d)$ \refpp{bw0}.  
The {\it standard generating set} is all 
 of vectors of the form 
$\frac 1{2^m}v_A$, where $m$ is a nonnegative integer and $A$ is an affine $2m$-space in $\Omega$.  In fact, this is just the set of minimal vectors of $\bw d$.  
\end{de}

\begin{prop}\labtt{minvecbw}
The minimal vectors in $\bw d$ 
are of the form $\frac 1{2^m}v_A \vep_S$,  where $m$ is a nonnegative integer, $0\le m \le 2^{\lfloor \frac d2 \rfloor}$,  $A$ is an affine $2m$-space in $\Omega$ and $S\in RM(2,d)$.  
They have norms 
$2^{\lfloor \frac d2 \rfloor}$.  
\end{prop} 
\pf
This is a standard result \cite{be,bwy}.  
\eop

\begin{de}\labttr{lowerframe}  Let $L:=\bw d$.  A {\it lower frame} or a {\it standard frame} is a set of $2^{d+1}$ minimal vectors of $L$ which forms an orbit under the action of the normal extraspecial subgroup 
of order $2^{1+2d}$ of $\brw d$.  (A lower frame was called a sultry frame in \cite{ibw1}.)   A {\it standard basis} or a {\it lower basis} is a basis contained in a standard frame with a labeling by $\Omega$ such that the set of minimal vectors of $L$ is as described in \refpp{minvecbw}.  An arbitrary labeling by $\Omega$ 
of a basis contained in a frame may not have this property.  See \cite{bwy}.  
\end{de}

\subsection{Review of commutator density} 

This concept was introduced in \cite{bwy}.  
Let $D$ be an extraspecial 2-group and let $Mod(D,-)$ be the category of modules for which the central involution of $D$ acts as $-1$.   
Often, $D$ is dihedral of order 8.  

The basic results are summarized in this section.  
For a proof, see \cite{bwy}. 

\begin{de}\labttr{commdens} 
Let $E$ be a group, $S$ a subset of $E$ and $M$ a $\ZZ [E]$ module.  We say that $S$ is {\it commutator dense on $M$} if $[M,E]=[M,S]$. 
\end{de}  

\begin{de}\labttr{2/4} 
Let $D$ be a dihedral group of order 8 and let $M$ be a $\ZZ [D]$-module.  We say that $M$ has the {\it 2/4 generation property} if for any pair of involutions $u, v$ which generate $D$, we have $M^+(u)+M^+(v)=M$.  
\end{de}

\begin{prop}\labtt{commdens2}
Let $D$ be a dihedral group of order 8 and let $M$ be a $\ZZ [D]$-module on which the central involution of $D$ acts as $-1$.  
Let $f\in D$ have order 4.  
Then on $M$, 
2/4-generation and commutator density of $\{ f\}$ are equivalent.  
\end{prop}
\pf \cite{bwy}.  
\eop

\begin{nota}\labttr{twist} 
Suppose that $D$ is dihedral of order 8 and that  $L$ is in the category $Mod(D,-)$.    Let $f$ be an element of order 4 in $D$ and let $p$ be an integer.  The {\it $p$-th twist} of $L$ is the $D$-submodule $L[p]:=L(f-1)^p$ of $\QQ \otimes L$. 
\end{nota}  

\begin{prop}\labtt{commdensapp}
Let $L=\bw d$ and let $f\in \rd d$ be a fourvolution.  Then $[L,\rd d]=L(f-1)$, i.e., $f$ is commutator dense on the $\rd d$-module $L$.  
\end{prop}
\pf \cite{bwy}.  \eop 

\begin{rem}\labttr{commdensapp2} 
The notation $L[p]$ (rather than $L(f-1)$) stresses dependence on $\rd d$ rather than on choice of fourvolution $f\in \rd d$ \refpp{commdensapp}.  This independence can be useful.  
\end{rem}

\section{Appendix: the minimal vectors of $\bw d [1]$}

The minimal vectors of 
$\bw d$ constitute the standard generating set \refpp{standardgenset}, as is well-known.  
We need the following fact about twists of Barnes-Wall lattices.  
This result may be new.

\begin{thm}\labtt{minvecbwd[1]}
The set of minimal vectors of $\bw d [1]$ is $K:=\cup_{m\ge 0}
K_m$, where $K_m$ is the set of all 
$2^{-m}v_A\vep _S$, where $A$ is a $(2m+1)$-dimensional affine subspace of $\Omega = \FF_2^d$ and $S \in RM(2,d)$.  
\end{thm}
\pf 
Define $L:=\bw d$. 
We use the commutator density property, that $L[1]$ equals 
$L(\pm f \pm 1)$
for any lower fourvolution $f$ \refpp{commdensapp}.   

Let $J$ be the set of minimal vectors in $L$.  Since each $f-1$ doubles norms and 
maps $L$ onto $L[1]$, it takes $J$ onto the set $K'$ of minimal vectors of $L[1]$.

The $K_m$ are orbits for the action of the standard monomial subgroup of $\brw d$.  To prove $K \subseteq K'$, it suffices to prove that $J(f-1)$ contains a single member of each $K_m$.  
It suffices to prove that, given $m$ such that $K_m \ne \emptyset$, that there exists a lower fourvolution $f$ so that $K_m \cap J(f-1)\ne \emptyset$.  

Take $A$, an affine $(2m+1)$-dimensional space.  Let $H$ be a hyperplane such that $dim(A\cap H)=2m$.  
Let $\tau$ be a translation on $\Omega$ which fixes $A$ and interchanges $H$ and $H+\Omega$.  
Define $f:=\tau \vep_H$, a lower fourvolution.  
Then $2^{-m}v_{A\cap H} \in J$ and 
$2^{-m}v_{A\cap H}(1+f)=2^{-m}v_A$.  

Finally, to prove that $K'\le K$, observe that if 
$v\in K$, the vector $v(f-1)^{-1} \in J$, so has the form 
$u=2^{-m}v_B\vep_S$, for some affine $2m$-space $B$.   
Then $u\tau \vep_H=(2^{-m}v_{B\tau} - 2^{-m+1}v_{B\tau \cap H})\vep_S$.  

If $B=B\tau$, $v=u(f-1)=2^{-m+1}v_{B\cap H}\vep_S\in K_{m-1}$.  

If $B\ne B\tau$, then $B\cap B\tau = \emptyset$  and 
$v=u(f-1)=
2^{-m}_{B+B\tau}\vep_{S+H} \in K_{m1}$.  
\eop


\end{document}